\documentclass [11pt,oneside,a4paper,mathscr]{amsart}

\usepackage{amscd,amsmath,amssymb,euscript,graphics}
\usepackage[frame,cmtip,curve,arrow,matrix,line,graph]{xy}

\title{Stability conditions on a non-compact Calabi-Yau threefold}
\author{Tom Bridgeland}

\date{}

\newtheorem{thm}{Theorem}[section]
\newtheorem{conj}[thm]{Conjecture}

\newtheorem{prop}[thm]{Proposition}
\newtheorem{lemma}[thm]{Lemma}
\newenvironment{pf}{\paragraph{Proof}}{\qed\par\medskip}
\theoremstyle{definition}
\newtheorem{defn}[thm]{Definition}

\newcommand{\conf}{\operatorname{C}}
\newcommand{\conft}{\operatorname{\tilde{C}}}
\newcommand{\Sym}{\operatorname{Sym}}
\newcommand{\U}{\mathcal{U}}

\newcommand{\Stabredo}{\operatorname{Stab}^{0}_{n}}
\newcommand{\Stabo}{\operatorname{Stab}^{0}}
\newcommand{\into}{\hookrightarrow}
\newcommand{\SL}{\operatorname{SL}}

\newcommand{\nab}{\check{\nabla}}

\renewcommand{\P}{\mathcal{P}}
\newcommand{\Stab}{\operatorname{Stab}}

\renewcommand{\leq}{\leqslant}
\renewcommand{\geq}{\geqslant}

\newcommand{\Lie}{\operatorname{Lie}}
\renewcommand{\Re}{\operatorname{Re}}
\renewcommand{\Im}{\operatorname{Im}}

\newcommand{\K}{\operatorname{K}}
\newcommand{\A}{\mathcal{A}}
\renewcommand{\AA}{\mathbb{A}}

\newcommand{\T}{\mathcal T}
\newcommand{\F}{\mathcal F}

\newcommand{\ZZ}{\mathcal{Z}}
\newcommand{\Aut}{\operatorname{Aut}}

\newcommand{\tensor}{\otimes}
\newcommand{\C}{\mathbb C}
\newcommand{\Z}{\mathbb Z}
\newcommand{\PP}{\mathbb P}

\newcommand{\id}{\operatorname{{1}}}
\newcommand{\D}{\operatorname{\mathcal{D}}}
\newcommand{\R}{\operatorname{R}}
\renewcommand{\L}{\operatorname{L}}
\newcommand{\WW}{\mathcal{W}}

\newcommand{\blob}{{\scriptscriptstyle{\bullet}}}

\newcommand{\Coh}{\operatorname{Coh}}
\newcommand{\OO}{\mathcal O}
\newcommand{\Mt}{\tilde{M}}

\newcommand{\End}{\operatorname{End}}

\newcommand{\Hom}{\operatorname{Hom}}

\newcommand{\eu}{\operatorname{\chi}}
\newcommand{\lra}{\longrightarrow}

\newcommand{\Mod}{\operatorname{Mod}}

\newcommand{\bigmat}[9]{\left(\begin{array}{ccc}#1 & #2 & #3 \\ #4 & #5 & #6 \\ #7 & #8 & #9 \end{array} \right)}

\newcommand{\RHom}{\operatorname{Hom}}
\newcommand{\ha}{\frac{1}{2}}
\newcommand{\mt}{\tilde{m}}

\begin{document}
\begin{abstract}
We study the space of stability conditions $\Stab(X)$ on the non-compact Calabi-Yau threefold $X$
which is the total space of the canonical bundle of $\PP^2$. We give a combinatorial description of an open subset of $\Stab(X)$ and state a conjecture
relating $\Stab(X)$ to the Frobenius manifold obtained from the quantum cohomology of $\PP^2$. We give some evidence from mirror symmetry for this conjecture. 
\end{abstract}

\maketitle


\section{Introduction}

The space of stability conditions $\Stab(X)$ on a variety $X$ was introduced in \cite{B} as a mathematical framework for understanding Douglas's notion of $\pi$-stability for D-branes in string theory \cite{Do4}.
This paper is concerned with the case when $X=\OO_{\PP^2}(-3)$ is the total space of the canonical line bundle of $\PP^2$. This non-compact Calabi-Yau threefold provides an amenable but interesting example on which to test the general theory, and many  features of the spectrum of D-branes on $X$ have already been studied in the physics literature (see for example \cite{DG,Do2,FHHI}).

So far, we are unable to give a complete description of  $\Stab(X)$. However, using the results of \cite{B2}, we define an open subset $\Stabo(X)\subset \Stab(X)$ which is a disjoint union of regions indexed by the elements of an affine braid group. 
The combinatorics of these regions leads us to conjecture a precise connection between $\Stab^0(X)$ and the Frobenius manifold defined by the quantum cohomology of $\PP^2$. Our main aim is to assemble some convincing evidence for this conjecture and to discuss some of its consequences.

The existence of deep connections between quantum cohomology and derived categories has been known for some time. In particular, following observations of Cecotti and Vafa \cite{CV} and Zaslow \cite{Za}, Dubrovin conjectured \cite{Du2} that the derived category of a Fano variety $Y$ has a full exceptional collection $(E_0,E_1,\cdots, E_{n-1})$ if and only if the quantum cohomology of $Y$ is generically semisimple, and that in this case the Stokes matrix $S_{ij}$ of the corresponding Frobenius manifold should coincide with the Gram matrix $\chi(E_i,E_j)$ for the Euler form of $\D(Y)$. 
This statement has been verified for projective spaces \cite{Gu,Ta}.

It was pointed out by Bondal and Kontsevich that a heuristic explanation for Dubrovin's conjecture can be given using mirror symmetry.
The mirror of a Fano variety with a full exceptional collection  is expected to be  an affine variety $\check{Y}$, together with a holomorphic function $f\colon \check{Y}\to \C$ with isolated singularities. The Frobenius manifold arising from the quantum cohomology of $Y$ should coincide with the Frobenius manifold of Saito type defined on the universal unfolding space of $f$. The Stokes matrix is then the intersection form evaluated on a distinguished basis of vanishing cycles. Under Kontsevich's homological mirror  proposal \cite{Ko} the intersection form is identified with the Euler form on $\D(Y)$, and the vanishing cycles, which are discs, correspond to exceptional objects in $\D(Y)$. 

The conjecture stated below suggests that it may be possible to use spaces of stability conditions to give
a more direct link between derived categories and quantum cohomology. To make this work one should somehow define the structure of a Frobenius manifold on the space of stability conditions which in some small patch recovers the usual quantum cohomology picture. At present however, the author has no clear ideas as to how this could be done.

The other general conclusion one can draw from the example studied in this paper is that the space of stability conditions $\Stab(X)$ is not an analogue of the stringy K\"{a}hler moduli space, but rather some extended version of it. The picture seems to be that the space $\Stab(X)$ is a global version of the  Frobenius manifold defined by big quantum cohomology, and the stringy K\"{a}hler moduli space is a submanifold which near the large volume limit is defined by the small quantum cohomology locus.

In the rest of the introduction we shall describe our results in more detail. Missing definitions and proofs are hopefully covered in the main body of the paper.
\subsection{}
A stability condition \cite{B} on a triangulated category $\D$ consists of a full abelian subcategory $\A\subset \D$ called the heart, together with a group homomorphism
\[Z\colon K(\D)\lra \C\]
called the central charge,
with the compatibility property that for every nonzero object $E\in \A$ one has
\[ Z(E)\in H=\{z\in \C : z=r\exp(i\pi \phi)\text{ with }r>0 \text{ and } 0<\phi\leq 1\}.\]
One insists further that $\A\subset \D$ is the heart of a bounded t-structure on $\D$, and that the map $Z$ has the Harder-Narasimhan property. 
The set of all stability conditions on $\D$ satisfying an extra condition called local-finiteness form a complex manifold $\Stab(\D)$.
Forgetting the heart $\A\subset \D$ and remembering the central charge gives a map
\[\ZZ\colon \Stab(\D) \lra \Hom_{\Z}(\K(\D),\C).\]
In this paper we shall consider the case when $\D$ is the subcategory of the bounded derived category of coherent sheaves on $X=\OO_{\PP^2}(-3)$ consisting of complexes whose cohomology sheaves are supported on the zero section $\PP^2\subset X$. In that case the Grothendieck group $\K(\D)$ is a free abelian group of rank three.
Our main result is

\begin{thm}
\label{one}
There is a connected open subset $\Stabo(X)\subset \Stab(X)$ which can be written as a disjoint union of regions
\[\Stabo(X)=\bigsqcup_{g\in G} D(g),\]
where $G$ is the affine braid group with presentation
\[G=\big\langle \tau_0,\tau_1,\tau_2 \;|\;\tau_i\tau_j\tau_i=\tau_j\tau_i\tau_j\text{ for all }i,j\big \rangle.\]
Each region $D(g)$ is mapped isomorphically by $\ZZ$ onto a locally-closed subset of the three dimensional vector space $\Hom_{\Z}(\K(\D),\C)$, and the closures of two regions $D(g_1)$ and $D(g_2)$ intersect in $\Stab^0(X)$ precisely if $g_1 g_2^{-1}=\tau_i^{\pm 1}$ for some $i\in \{0,1,2\}$.
\end{thm}

The stability conditions in a given region $D(g)$ all have the same heart $\A(g)\subset \D$. Each of these categories $\A(g)$ is equivalent to a category of nilpotent representations of a quiver with relations of the form
\[\includegraphics{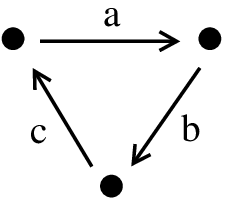}\]
where the positive integers $a,b,c$ labelling the graph represent the number of arrows in the quiver joining the corresponding vertices. In fact the triples $(a,b,c)$ which come up are precisely the positive integer solutions to the Markov equation
\[a^2 + b^2 + c^2 =abc.\]
We denote by $S_0(g),S_1(g),S_2(g)$ the three simple objects of $\A(g)$ corresponding to the three one-dimensional representations of the quiver. In the case when $g=e$ is the identity we simply write $S_i=S_i(e)$. The objects $S_i(g)$ are spherical objects of $\D$ in the sense of Seidel and Thomas \cite{ST}. As such they define autoequivalences
$\Phi_{S_i(g)}\in \Aut \D$. These descend to give automorphisms
\[\phi_{S_i(g)}\in \Aut \K(\D)\quad i=0,1,2,\]
which with respect to the fixed basis of $\K(\D)$ defined by the classes of the objects $S_i$ are given by a triple of matrices \[P_0(g),P_1(g),P_2(g)\in \SL(3,\Z).\]
 It turns out that exactly the same system of matrices come up in the study of the quantum cohomology of $\PP^2$.

\subsection{}
Dubrovin showed that the semisimple Frobenius structure arising from the quantum cohomology of $\PP^2$ can be analytically continued to give a Frobenius structure on a dense open subset $M$ of the universal cover of the configuration space
\[\conf_3(\C)=\{(u_0,u_1,u_2)\in \C: i\neq j \implies u_i \neq u_j\}/\Sym_3.\]
Note that in some small ball on $M$ the corresponding prepotential encodes the geometric data of the Gromov-Witten invariants of $\PP^2$, but away from this patch there is no such direct interpretation. Thus, just like the space of stability conditions, $M$ is a non-perturbative object, not depending on any choice of large volume limit.

Given a point $m\in M$ we denote by $\{u_0(m),u_1(m),u_2(m)\}$ the corresponding unordered triple of points in $\C$, and set
\[\C_m=\C\setminus\{u_0(m),u_1(m),u_2(m)\}.\]
Let $W$ denote the space
\[W=\{(m,z)\in M\times \C : z\in \C_m\}\]
with its projection  $p\colon W \to M$.
Using the Frobenius structure Dubrovin defined a series of flat, holomorphic connections $\check{\nabla}^{(s)}$ on the pullback of the tangent bundle $p^* (\T_M)$.
These connections are called the second structure connections. Connections of this type were first introduced by K. Saito in the theory of primitive forms for unfolding spaces. We shall be interested only in the case $s=\ha$.

For each $m\in M$ the connection $\nab=\nab^{(\ha)}$ restricts to give a holomorphic connection $\nab_m$ on a trivial rank three bundle over $\C_m$.
Dubrovin showed that this family of connections is isomonodromic.
Define another configuration space
\[\conf_3(\C^*)=\{(u_0,u_1,u_2)\in \C^*: i\neq j \implies u_i \neq u_j\}/\Sym_3\]
and let $\conft_3(\C^*)$ be its universal cover.
Define
\[M^0=\{m\in M:0\in \C_m\}\]
and let $\Mt^0$ be its inverse image in $\conft_3(\C^*)$ under the natural map $\conft_3(\C^*)\to\conft_3(\C)$. We can choose a base-point $m\in M^0$ such that  $\{u_0(m),u_1(m),u_2(m)\}$ are the three roots of unity. Let $(\gamma_0,\gamma_1,\gamma_2)$ denote the following basis of $\pi_1(\C_m,0)$.
\[\includegraphics{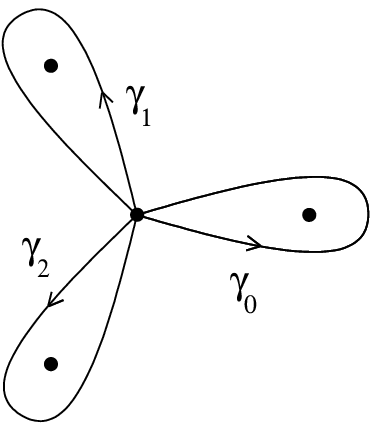}\]
Let $m\in U\subset M^0$ be a small simply-connected neighbourhood of $m$. For each point $m'\in U$ there is a chosen basis of $\pi_1(\C_{m'},0)$ obtained  by deforming the loops $\gamma_i$, which we also denote $(\gamma_0,\gamma_1,\gamma_2)$.

Let $V$ be the space of flat sections of $\nab_m$ near the origin $0\in \C$. 
Using the connection $\nab$ we can identify $V$ with the space of flat sections of $\nab_{\mt}$ near $0\in \C$ for all points $\mt\in \Mt^0$.
As we explain in Section \ref{braid}, the group $G$ is a subgroup of $\pi_1(\conf_3(\C^*))$, and hence acts on $\conft_3(\C^*)$. Taking the monodromy of the connection $\nab_{\mt}$ around the loops $\gamma_i$ for $\mt\in g(U)\cap\Mt^0$ we obtain linear automorphisms
\[\alpha_i(g)\in \Aut(V)\]
for $i=0,1,2$. The following result relates these to the transformations $\phi_{S_i(g)}$ of the last section.

\begin{thm}
\label{two}
There is a triple of flat sections $(\phi_0,\phi_1,\phi_2)$ of the second structure connection $\nab$ such that for all $g\in G$ the monodromy transformations $\alpha_i(g)$ act by the matrices $P_i(g)$ defined in the last section. This condition fixes the triple $(\phi_0,\phi_1,\phi_2)$ uniquely up to a scalar multiple.
\end{thm}

This Theorem is a simple recasting of some results of Dubrovin, and boils down to two previously observed  coincidences. The first is the fact mentioned in the introduction that the Stokes matrix $S_{ij}$ for the quantum cohomology of $\PP^2$ coincides with the Gram matrix $\eu(E_i,E_j)$ for the Euler form on $\K(\PP^2)$ with respect to a basis coming from an exceptional triple of vector bundles $(E_0,E_1,E_2)$. The second is that this coincidence is compatible with the braid group actions on these matrices arising on the one hand  from the analytic continuation of the Frobenius manifold \cite[Theorem 4.6]{Du3}, and on the other from the action of mutations on exceptional triples discovered by Bondal, Gorodentsev and Rudakov \cite{Bo,GR}.

In fact the connection $\nab$ corresponds to the Gauss-Manin connection on the universal unfolding space of the singularity mirror to the space $X$. It might perhaps be easier to understand the connection in this geometric way. But part of the point of this paper is to try to avoid passing to the mirror.

\subsection{}
We now describe a conjecture relating the quantum cohomology of $\PP^2$ to the space of stability conditions on $X$. The noncompactness of $X$ makes this relationship slightly more complicated than might be expected. In particular, the Euler form 
$\chi(-,-)$ on $\K(\D)$ is degenerate, with a one-dimensional kernel generated by the class of a skyscraper sheaf $[\OO_x]\in \K(\D)$  for $x\in \PP^2\subset X$. In terms of the basis defined by the spherical objects $S_i=S_i(e)$ one has
\[[\OO_x]=[S_0]+[S_1]+[S_2].\]
Since this class is somehow special, and in particular is preserved by all autoequivalences of $\D$, it makes sense to 
define a space of normalised  stability conditions by
\[\Stabredo(X)=\{\sigma=(Z,\P)\in \Stabo(X): Z(\OO_x)=i\}.\]
This is a connected submanifold of $\Stabo(X)$. Define an affine space
\[\AA^2=\{(z_0,z_1,z_2)\in \C^3: z_0+z_1+z_2=i\}.\]
In co-ordinate form the map $\ZZ$ gives a local isomorphism
\[\ZZ\colon \Stabredo(X) \lra \AA^2\]
obtained by sending a stability condition to the triple $(Z(S_0),Z(S_1),Z(S_2))$.

On the quantum cohomology side, the flat sections $(\phi_0,\phi_1,\phi_2)$ of Theorem \ref{two} do not form a basis, and in fact satisfy $\phi_0+\phi_1+\phi_2=0$. Pulling back the connection $\nab$ via the embedding
$M^0\to W$ defined by $p\mapsto (p,0)$
we obtain a flat connection on the tangent bundle $\T_{M^0}$. Taking co-ordinates whose gradients are the sections $(\phi_0,\phi_1,\phi_2)$ and rescaling appropriately, one obtains a holomorphic map
\[\WW\colon \Mt^0 \lra \AA^2.\]
This map is invariant under the free $\C$ action on $\Mt^0\subset \conft_3(\C^*)$ which lifts the $\C^*$ action on $\conf_3(\C^*)$ obtained by simultaneously rescaling the points $(u_0,u_1,u_2)$.
The quotient $\conft_3(\C^*)/\C$ is the universal covering space of
\[\{[u_0,u_1,u_2]\in \PP^2: u_i\neq 0\text{ and }u_i\neq u_j\}.\]
The quotient $\Mt^0/\C$ is therefore a dense open subset. 
We call the induced map
\[\WW\colon \Mt^0/\C \lra \AA^2.\]
the homogeneous twisted period map. It is a local isomorphism. We can now state our conjecture.
 
\begin{conj}
\label{conj}
There is a commuting diagram
\[\begin{CD}
\Stabredo(X) &@>F>> &\Mt^0/\C\\
@V{\ZZ}VV && @V{\WW}VV \\
\AA^2 &@= &\AA^2
\end{CD}\]
Moreover $F$ is an isomorphism onto a dense open subset.
\end{conj}

Proving this conjecture would require a more detailed understanding of the geometry of the homogeneous twisted period map. In particular, it would be necessary to find an open subset of $\Mt^0/\C$ which was mapped isomorphically by $\WW$ onto the subset
\[\{(z_0,z_1,z_2)\in \AA^2: \Im(z_i)>0\}\]
which is the image of the interior of the region $D(e)$ under the map $\ZZ$.

\subsection{}
Here we describe two pieces of evidence for Conjecture \ref{conj}. First
consider the submanifold $D\subset \conf_3(\C)$ defined parametrically by taking the unordered triple of points
\[u_i=-1+z^{1/3}\]
for some $z\in \C\setminus\{0,1\}$. The inverse image of $D$ in the universal cover $\conft_3(\C)$ is contained in the open subspace $M$. The submanifold $D$ (or its inverse image in $M$) is the small quantum cohomology locus; in the standard flat co-ordinates  it is given by $(t_0,t_1,t_2)=(-1,e^z,0)$. Dubrovin showed \cite[Proposition 5.13]{Du5} that on this locus the homogeneous twisted period map satisfies the differential equation
\[\bigg[\theta_z^3-z(\theta_z+\frac{1}{3})(\theta_z+\frac{2}{3})\theta_z\bigg]\WW=0
\qquad\theta_z\equiv z\frac{d}{dz}.\]
This is the Picard-Fuchs equation for the periods of  the mirror of $X$, and is thus precisely the equation satisfied by the central charge on the stringy K\"ahler moduli space \cite{AGM,DG}.

A second piece of evidence for Conjecture \ref{conj} is that if we go down a dimension to the case $X=\OO_{\PP^1}(-2)$ the corresponding statement is known to be at least nearly true.
In that case the space $M$
is the universal cover of
\[\conf_2(\C)=\{(u_0,u_1)\in \C:u_0\neq u_1\}\]
so that $\Mt^0/\C$ is the universal cover of $\C\setminus\{0,1\}$ with co-ordinate $\lambda=u_1/u_0$. In this case we must take the second structure connection with $s=0$ (in general, for a projective space of dimension $d$ we should take $s=(d-1)/2$). Thus the homogeneous twisted period map in this case is just the homogeneous part of the standard period map for the quantum cohomology of $\PP^1$. This was computed by Dubrovin. Identifying the affine space
\[\AA^1=\{(z_0,z_1):z_0+z_1=i\}\]
with $\C$ via the map $(z_0,z_1)\mapsto z_0$, the equation \cite[G.20]{Du1} implies that the homogeneous period map is
\[\WW(\lambda)=\bigg(\frac{1}{\pi}\bigg)\cos^{-1}\bigg(\frac{1+\lambda}{1-\lambda}\bigg).\]

On the other hand the space $\Stab(X)$ was studied in \cite{Br4}. The corresponding open subset $\Stabo(X)\subset\Stab(X)$  is actually  a connected component, and the corresonding space $\Stabredo(X)$  is a covering space of $\C\setminus \Z$. This gives the following result.

\begin{thm}
In the case $X=\OO_{\PP^1}(-2)$
there is a commuting diagram
\[\begin{CD}
\Stabredo(X) &@<H<< &\Tilde{\C\setminus\{0,1\}}\\
@V{\ZZ}VV && @V{\WW}VV \\
\C\setminus \Z &@= &\C\setminus \Z
\end{CD}\]
in which all the maps are covering maps.
\end{thm}

In fact one expects $\Stabredo(X)$ to be simply-connected so that $H$ is actually an isomorphism.

\subsection*{Acknowledgements}
The problem of describing $\Stab(\OO_{\PP^2}(-3))$ was originally conceived as a joint project with Alastair King, and the basic picture described in Theorem \ref{one} was worked out jointly with him. It's a pleasure to thank Phil Boalch who first got me interested in the connections with Stokes matrices and quantum cohomology. Several other people have been extremely helpful in explaining various things about Frobenius manifolds; let me thank here B. Dubrovin, C. Hertling and M. Mazzocco. 


\section{Stability conditions on $X$}

In this section we  justify the claims about $\Stab(X)$ made in  the introduction. In particular we  prove Theorem \ref{one}. We start by summarising some of the necessary definitions. More details  can be found in \cite{B,B2}.

\subsection{Stability conditions and tilting}

Let $\D$ be a triangulated category.
Recall that a bounded t-structure on $\D$ determines and is determined by its heart, which is an abelian subcategory $\A\subset \D$.
One has an identification of Grothendieck groups $\K(\D)=\K(\A)$.

A stability function on an abelian category $\A$ is defined to be a  group homomorphism $Z\colon \K(\A)\to\C$
such that
\[0\neq E\in\A \implies Z(E)\in\mathbb{R}_{>0}\,\exp({i\pi\phi(E)})\text{ with
}0<\phi(E)\leq 1.\]
The real number $\phi(E)\in(0,1]$ is called the phase of the object $E$.

A nonzero object $E\in\A$ is said to be
semistable
with respect to $Z$
if every subobject $0\neq A\subset E$ satisfies $\phi(A)\leq\phi(E)$.
The stability function $Z$ is said to have the Harder-Narasimhan property if
every nonzero object $E\in\A$ has
a finite filtration
\[0=E_0\subset E_1\subset \cdots\subset E_{n-1}\subset E_n=E\]
whose factors $F_j=E_j/E_{j-1}$ are semistable objects of $\A$ with
\[\phi(F_1)>\phi(F_2)>\cdots>\phi(F_n).\]
A simple sufficient condition for the existence of Harder-Narasimhan
filtrations was given in \cite[Proposition 2.4]{B}.
In particular the Harder-Narasimhan
property
always holds when $\A$ has finite length.

The definition of a stability condition appears in \cite{B}. For our purposes the following equivalent definition will be more useful,
see \cite[Proposition 5.3]{B}.

\begin{defn}
A stability condition on  $\D$ 
consists of a bounded t-structure on $\D$ and a stability function on its heart which
has the Harder-Narasimhan property.
\end{defn}

The induced map $Z\colon \K(\D)\to\C$ is called the central charge of the stability condition.
It was shown in \cite{B} that the set of stability conditions on $\D$ satisfying an additional condition called local-finiteness form the points of a complex manifold $\Stab(\D)$.
In general this manifold will be infinite-dimensional, but in the cases we  consider in this paper $\K(\D)$ is of finite rank, and it follows that $\Stab(\D)$ has finite dimension.

To construct t-structures we use the method of tilting introduced by Happel, Reiten and Smal\o \,\cite{HRS}, based on earlier work of Brenner and Butler \cite{BB}.
Suppose $\A\subset \D$ is the heart of a bounded t-structure and is a
finite length abelian category. Note that the t-structure is completely determined by the set of simple objects of $\A$; indeed $\A$ is the smallest extension-closed subcategory of $\D$ containing this set of objects. Given a simple object $S\in \A$ define
$\langle S \rangle\subset \A$ to be the full subcategory consisting of objects $E\in\A$ all
of whose simple factors are isomorphic to $S$. One can either view $\langle S \rangle$ as
the torsion part of a torsion theory on $\A$, in which case the torsion-free
part is
\[\F=\{E\in \A:\Hom_{\A}(S,E)=0\},\]
or as the torsion-free part, in which case the torsion part is
\[\T=\{E\in\A:\Hom_{\A}(E,S)=0\}.\]
The corresponding tilted subcategories are defined to be
\begin{eqnarray*}
\L_S \A &=& \{E\in\D:H^i(E)=0\text{ for
}i\notin\{0,1\},H^{0}(E)\in\F\text{ and }H^1(E)\in\langle S \rangle\}  \\
\R_S \A &=& \{E\in \D:H^i(E)=0\text{ for
}i\notin\{-1,0\},H^{-1}(E)\in\langle S \rangle\text{ and }H^0(E)\in\T\}.
\end{eqnarray*}
They are the hearts of new bounded t-structures on $\D$.

\subsection{Quivery subcategories}
Let $X=\OO_{\PP^2}(-3)$ with its projection $\pi\colon X\to \PP^2$.
Let $\D$ denote the full subcategory of the bounded derived category
of coherent sheaves on $X$ consisting of complexes supported on the zero-section $\PP^2\subset X$.
Let $\Stab(X)$ denote the space of locally-finite stability conditions on $\D$.

Let $(E_0,E_1,E_2)$ be an exceptional collection of vector bundles on $\PP^2$. Any exceptional collection in $\D^b\Coh(\PP^2)$ is of this form up to shifts. It was proved in \cite{B2} that there is an equivalence of categories
\[ \RHom^{\blob}\big(\bigoplus_{i=0}^{2} \pi^* E_i, -\big)\colon
\D^b\Coh(X) \lra \D^b\Mod (B),\]
where $\Mod(B)$ is the category of finite-dimensional right modules for the algebra
\[B=\End_X\big(\bigoplus_{i=0}^{2} \pi^* E_i\big).\]
The algebra $B$ can be described as the path algebra of a quiver with relations taking the form
\[\includegraphics{quiver.eps}\]
Pulling back the standard t-structure on $\D^b \Mod(B)$ gives a bounded t-structure on $\D$
whose heart is equivalent to the category of nilpotent modules of $B$. The abelian subcategories $\A\subset \D$ obtained in this way are called exceptional. An abelian subcategory of $\D$ is called quivery if it is of the form $\Phi(\A)$ for some exceptional subcategory $\A\subset \D$ and some autoequivalence $\Phi\in \Aut (\D)$.

Any quivery subcategory $\A\subset \D$ is equivalent to a category of nilpotent modules of an algebra of the above form. As such it has three simple objects $\{S_0,S_1,S_2\}$ corresponding to the three one-dimensional representations of the  quiver. These objects $S_i$ are spherical in the sense of Seidel and Thomas and thus give rise to autoequivalences
\[\Phi_{S_i}\in \Aut(\D).\]
Note that the three simples $S_i$ completely determine the corresponding quivery subcategory $\A\subset \D$. The Ext groups between them can be read off from the quiver
\[\Hom^1_{\D}(S_0,S_1)=\C^a, \quad \Hom^1_{\D}(S_1,S_2)=\C^b, \quad \Hom^1_{\D}(S_2,S_0)=\C^c\]
with the other $\Hom^1$ groups being zero. Serre duality then determines the other groups.

Take $\A$  to be the exceptional subcategory of $\D$ corresponding to the exceptional collection $(\OO,\OO(1),\OO(2))$ on $\PP^2$. Its simples are
\[S_0=i_* \OO, \quad S_1=i_* \Omega^1(1)[1], \quad S_2=i_* \OO(-1)[2],\]
where $i\colon \PP^2\into X$ is the inclusion of the zero-section, and $\Omega$ denotes the cotangent bundle of $\PP^2$. We have $(a,b,c)=(3,3,3)$.

Let us compute the automorphisms $\phi_{S_i}$ of $\K(\D)$ induced by the autoequivalences $\Phi_{S_i}$. The twist functor $\Phi_S$ is defined
by the triangle
\[\Hom_{\D}^{\blob}(S,E)\tensor S_i\lra E \lra \Phi_{S}(E)\]
so that, at the level of K-theory,
\[\phi_S([E])=[E]-\chi(S,E)[S].\]
If we write $P_i$ for the matrix representing the transformation $\phi_{S_i}$ with respect to the basis $([S_0],[S_1],[S_2])$ of $\K(\D)$ then
\[P_0=\bigmat{1}{3}{-3}{0}{1}{0}{0}{0}{1},\quad P_1=\bigmat{1}{0}{0}{-3}{1}{3}{0}{0}{1},
\quad P_2=\bigmat{1}{0}{0}{0}{1}{0}{3}{-3}{1}.\]

\subsection{Braid group action}
\label{braid}

It was shown in \cite{B2} that if one tilts a quivery subcategory $\A\subset \D$ at one of its simples one obtains another quivery subcategory. To describe this process in more detail we need to define a certain braid group which acts on triples of spherical objects.

The three-string annular braid group $CB_3$ is the fundamental group of the configuration space of three unordered points in $\C^*$.
It is generated by three elements $\tau_i$ indexed by the cyclic group $i\in \Z_3$ together with a single element $r$, subject to the relations
\begin{eqnarray*}
r\tau_i r^{-1}&=&\tau_{i+1} \text{ for all }i\in \Z_3, \\
\tau_i\tau_{j}\tau_i&=&\tau_{j}\tau_i\tau_{j}\text{ for all }i,j\in
\Z_3.
\end{eqnarray*}
For a proof of the validity of this presentation see \cite{KP}.
If we take the base point to be defined by the three roots of unity, then the elements $\tau_1$ and $r$ correspond to the
loops obtained by moving the points as follows:
\bigskip
\[\includegraphics{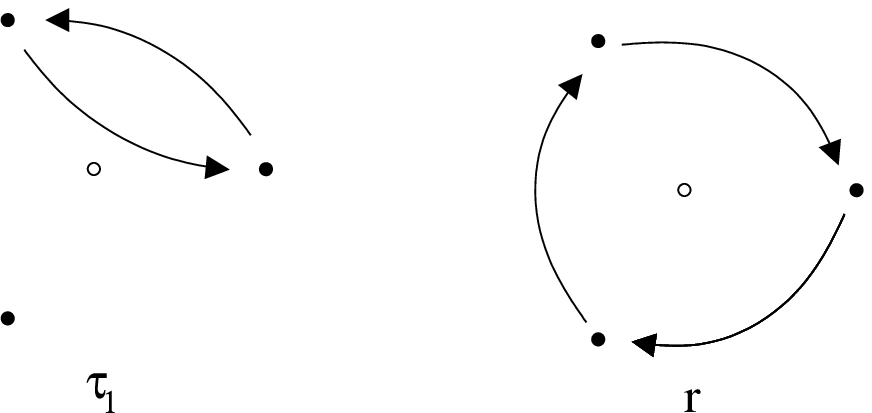}\]
We write $G\subset CB_3$ for the subgroup generated by the three braids $\tau_0, \tau_1,\tau_2$.

Define a spherical triple in  $\D$ to be a triple of spherical objects $(S_0,S_1,S_2)$ of $\D$. 
The group $CB_3$ acts on the set of spherical triples in $\D$ by the formulae
\[\tau_1(S_0,S_1,S_2)=(S_1[-1],\Phi_{S_1}(S_0), S_2),\quad r(S_0,S_1,S_2)=(S_2,S_0,S_1).\]
The following result allows one to completely understand the process of tilting for quivery subcategories of $\D$.

\begin{prop}
\label{tom}
Let $\A\subset \D$ be a quivery subcategory with simples $(S_0,S_1,S_2)$. Then for each $i=0,1,2$ the three 
simples of the tilted quivery subcategory  $\L_{S_i}(\A)$ are given by the spherical triple $\tau_i (S_0,S_1,S_2)$.
\end{prop}

For each $g\in G$ we then have a quivery subcategory $\A(g)\subset \D$ obtained by repeatedly tilting starting at $\A$. Its three simples are given by the spherical triple
\[(S_0(g), S_1(g), S_2(g))=g(S_0,S_1,S_2).\]
Note that the three simples of an arbitrary quivery subcategory have no well-defined ordering, but the above definition gives a chosen order for the simples of the quivery subcategories $\A(g)$.

Let $P_i(g)\in \SL(3,\Z)$ be the matrix representing the automorphism of $\K(\D)$ induced by the twist functor $\Phi_{S_i(g)}$ with respect to the fixed basis $([S_0],[S_1],[S_2])$. The formulae defining the action of the braid group on spherical triples show that this system of matrices has the following transformation laws
\begin{eqnarray*}
P_0(\tau_1 g)=P_1(g), \quad &P_1(\tau_1 g)=P_1(g) P_0(g) P_1(g)^{-1}, \quad &P_2(\tau_1 g)=P_2(g), \\
P_0(rg)=P_2(g), \quad &P_1(rg)=P_0(g), \quad &P_2(rg)=P_0(g).
\end{eqnarray*}

Introduce a graph $\Gamma(\D)$ whose vertices are the quivery subcategories of $\D$, and in which two subcategories are joined by an edge if they differ by a tilt at a simple object. It was shown in \cite{B2} that distinct elements $g\in G$ define distinct subcategories $\A(g)\subset \D$. It follows that each connected component of $\Gamma$ is just the Cayley graph of $G$ with respect to the generators $\tau_0,\tau_1,\tau_2$.

\subsection{Stability conditions on $X$}

Given an element $g\in G$ let $\A(g)\subset \D$ be the corresponding quivery subcategory.
The class of any nonzero object $E\in \A(g)$ is a strictly positive linear combination:
\[[E]= \sum n_i [S_i(g)] \text{ with }n_1,n_2,n_3 \geq 0\text{ not all zero}.\]
It follows that to define a stability condition on $\D$ we can just choose three complex numbers $z_i$ in the strict upper half-plane
\[H=\{z\in \C: z=r\exp(i\pi\phi) \text{ with }r>0\text{ and }0<\phi\leq 1\}\]
 and set $Z(S_i(g))=z_i$.
The Harder-Narasimhan property is automatically satisfied because $\A(g)$ has finite length.
We shall denote the corresponding stability condition by $\sigma(g,z_0,z_1,z_2)$.  

\begin{lemma}
If $\sigma=\sigma(g,z_0,z_1,z_2)$ is a stability condition on $\D$ of the sort defined above, and $E\in \D$ is stable in $\sigma$, then there is an open subset $U\subset \Stab(\D)$ containing $\sigma$ such that $E$ is stable for all stability conditions in $U$.
\end{lemma}

\begin{pf}
This follows from the arguments of \cite[Section 8]{Br}. It is enough to check that the set of classes $\gamma\in \K(\D)$ such that there is an object $F\in \D$ with class $[F]=\gamma$ such that $m_{\sigma}(F)\leq m_{\sigma}(E)$ is finite. This is easy to see because the heart of $\sigma$ has finite length.
\end{pf}

To each element $g\in G$ there is an associated set of stability conditions
\[D(g)=\{\sigma(g,z_0,z_1,z_2):(z_0,z_1,z_2)\in H^3 \text{ with at most one }z_i \in \mathbb{R}\}\subset \Stab(X).\]
By definition these subsets of $\Stab(X)$ are disjoint since they correspond to stability conditions with different hearts.

\begin{prop}
There is an open subset
\[\Stab^0(X)=\bigsqcup_{g\in G} D(g)\subset\Stab(X).\]
If $g_1,g_2\in G$ then the closures of the regions $D(g_i)$  intersect in $\Stab^0(X)$ precisely if $g_1=\tau_i^{\pm 1} g_2$ for some $i\in\{0,1,2\}$.
\end{prop}

\begin{pf}
Suppose a point $\sigma=\sigma(g,z_0,z_1,z_2)$ lies in $D(g)$. We must show that there is an open neighbourhood of $\sigma$ contained in the subset $\Stab^0(X)$.
The simple objects $S_i=S_i(g)\in \A(g)$ are stable in $\sigma$. They remain stable in a small open neighbourhood $U$ of $\sigma$ in $\Stab(X)$. We repeatedly use the easily proved fact that if $\A,\A'\subset\D$ are hearts of bounded t-structures and $\A\subset \A'$ then $\A=\A'$.

Suppose first that $\Im (z_i)>0$ for each $i$. Shrinking $U$ we can assume each $S_i$ has phase in the interval $(0,1)$ for all stability conditions $(Z,\P)$ of $U$. Since $\A(g)$ is the smallest extension-closed subcategory of $\D$ containing the $S_i$ it follows that $\A(g)$ is contained in the heart $\P((0,1])$ of all stability conditions in $U$. This implies that $\P((0,1])=\A(g)$ and so $U$ is contained in $D(g)$.

Suppose now that one of the $z_i$, without loss of generality $z_0$, lies on the real axis, so that $\sigma$ lies on the boundary of $D(g)$. Thus $z_0\in \mathbb{R}_{<0}$, and  $\Im(z_i)>0$ for $i=1,2$. Shrinking $U$ we can assume that $\Re Z(S_0)<0$ and $\Im Z(S_i)>0$ for $i=1,2$ for all stability conditions $(Z,\P)$ of $U$.

The object $S'=\Phi_{S_0}(S_2)\in \D$ lies in $\A(g)$, and is in fact a universal extension
\[0\lra S_2\lra S' \lra S_0^{\oplus a}\lra 0\]
where $a=\dim \Hom_{\D}^1(S_0,S_2)$.
Since $\Hom_{\D}(S_0,S')=0$ the object $S'$ lies in $\P((0,1))$ and shrinking $U$ we can assume that this is the case for all stability conditions $(Z,\P)$ of $U$.  

We split $U$ into the two pieces $U_+=\Im Z(S_0) \geq 0$ and $U_-=\Im Z(S_0)<0$. The argument above shows that $U_+\subset D(g)$. On the other hand, for any stability condition $(Z,\P)$ in $U_-$ the object $S_0$ is stable with phase in the interval $(1,3/2)$. Thus the heart $\P((0,1])$ contains the objects $S_0[-1], S'$ and $S_1$. Since these are the simples of the finite length category $\A(\tau_0 g)$ it follows that $U_-\subset D(\tau_0 g)$. 
\end{pf}


\section{Quantum cohomology and the period map}

In this section we describe some of Dubrovin's results concerning the twisted period map of the quantum cohomology of $\PP^2$.

\subsection{Frobenius manifolds}

The notion of a Frobenius manifold was first introduced by Dubrovin, although similar structures arising in singularity theory were studied earlier by K. Saito. A Frobenius manifold is a complex manifold $M$ with a flat metric $g$ and a commutative multiplication
\[\circ:\T_M\tensor \T_M\lra \T_M\]
 on its tangent bundle, satisfying the compatibility condition
\[g(X\circ Y,Z)=g(X,Y\circ Z).\]
One requires that locally on $M$ there exists a holomorphic function $\Phi$ called the prepotential such that
\[g(X\circ Y,Z)=XYZ(\Phi)\]
for all flat vector fields $X,Y,Z$.
Finally, one also assumes the existence of a flat identity vector field $e$, and an Euler vector field $E$ satisfying
\[\Lie_E(\circ) =\circ, \qquad \Lie_E(g)=(2-d) \cdot g\]
for some constant $d$ called the charge of the Frobenius manifold.

Given a smooth projective variety $Y$ of dimension $d$ one can put the structure of a Frobenius manifold of charge $d$ on an open subset of the vector space $H^*(Y,\C)$. The metric is the constant metric given by the Poincar\'{e} pairing, and the prepotential $\Phi$ is defined by an infinite series whose coefficients are the genus zero Gromov-Witten invariants, which naively speaking count rational curves in $Y$.
The condition that the resulting algebras be associative translates into the statement that $\Phi$ satisfies the WDVV equations. In turn, these equations boil down to certain relations between Gromov-Witten classes arising from the structure of the cohomology ring of the moduli space of pointed rational curves. The Gromov-Witten invariants are only non-vanishing in certain degrees, which gives the existence of the Euler vector field. The resulting Frobenius manifold is called the quantum cohomology of $Y$. 

Actually, for a general projective variety $Y$ (even a Calabi-Yau threefold) it is not known whether the series defining the prepotential has nonzero radius of convergence, so one has to work over a formal coefficient
ring.  But we shall only be interested in the case when $Y$ is a projective space and here it is known that there are no convergence problems.

For example, in the case $Y=\PP^2$ we can take co-ordinates $t=t_0 + t_1 \omega + t_2 \omega^2$ where $\omega\in H^2(X,\C)$ is the class of a line, and the function $\Phi$ is then defined by the series
\[\Phi(t)=\frac{1}{2}(t_0^2 t_2 + t_0 t_1^2) + \big.\sum_{k\geq 1} \frac{n_k}{(3k-1)!} \, t_2^{3k-1} e^{kt_1},\]
where $n_k$ is the number of curves of degree $k$ on $\PP^2$ passing through $3k-1$ generic points. The Euler vector field is
\[E=t_0 \frac{\partial}{\partial t_0} + 3\frac{\partial}{\partial t_1} - t_2\frac{\partial}{\partial t_2}.\]
See \cite[Lecture 1]{Du3} for more details.

\subsection{Tame Frobenius manifolds}

Let $M$ be a Frobenius manifold. Multiplication by the Euler vector field $E$ defines a section $\U\in \End(\T_M)$. A  point $m\in M$ is called tame if the endomorphism $\U$ has distinct eigenvalues. The set of tame points of $M$ forms an open (possibly empty) subset of $M$. A Frobenius manifold will be called tame if all its points are tame. 
Let
\[\conf_n(\C)=\{(u_0,\cdots,u_{n-1})\in \C^n: i\neq j \implies u_i\neq u_j\}/\Sym_n\]
be the configuration space of $n$ unordered points in $\C$.
Dubrovin showed that if $M$ is a tame Frobenius manifold the map $M \to \conf_n(\C)$ defined by the eigenvalues of $\U$ is a regular covering of an open subset of $\conf_n(\C)$.
This means that locally one can use the functions $u_i$ as co-ordinates on $M$. In terms of these canonical co-ordinates the product structure is
\[\frac{\partial}{\partial u_i} \cdot \frac{\partial}{\partial u_j}=\delta _{ij} \frac{\partial}{\partial u_i}.\]
and the Euler field takes the simple form
\[E=\big.\sum_i u_i \frac{\partial}{\partial u_i}.\]
The non-trivial data of the Frobenius structure on $M$ is entirely contained in the dependence of the metric on the canonical co-ordinates.

Given a Frobenius manifold $M$ it is natural to ask whether it is possible to analytically continue the prepotential $\Phi$ to obtain a larger Frobenius manifold $M'$ such that $M$ can be identified with an open subset of $M'$. Dubrovin showed how to do this for tame Frobenius manifolds.

\begin{thm}\cite[Theorem 4.7]{Du3}
\label{dub}
Given a tame Frobenius manifold $M$ of dimension $n$, there is a regular covering space $\conft_n(\C)\to \conf_n(\C)$, a divisor $B\subset \conft_n(\C)$, and a tame Frobenius structure on $M'=\conft_n(\C)\setminus B$ such that there is an open inclusion of Frobenius manifolds $M\into M'$.
\qed
\end{thm}

Let $M$ be a Frobenius manifold. Define a subset of $M\times \C$
\[W=\{(p,z)\in M\times \C:\det(\U-z\id)\neq 0 \},\]
and let $p\colon W \to M$ be the projection map. For each parameter $s\in \C$ one can define a flat, holomorphic connection $\nab=\nab^{(s)}$ on the bundle $p^* (T_M)$, by the following formulae
\begin{eqnarray*}
\nab_X^ Y&=&\nabla_X Y - (\nabla E + c\id)(U-z\id)^{-1}(X\circ Y) \\
\nab_{\partial/\partial z} Y&=&\nabla_{\partial/\partial z} Y + (\nabla E +c\id)(U-z\id)^{-1}(Y)
\end{eqnarray*}
Here $\nabla$ is the Levi-Civita connection corresponding to the flat metric $g$ on $M$, $\nabla E$ is the endomorphism of $\T_M$ defined by $X \mapsto \nabla_X E$, and
\[c= s+ \frac{(d-1)}{2}.\]
We shall be interested in the case when $c=d-1$.

Assume now that $M$ is a tame Frobenius manifold with its canonical co-ordinates $u_0,\cdots, u_{n-1}$. The space $W$ takes the form
\[W=\{(m,z)\in M\times \C : z\neq u_i\text{ for }0\leq i\leq n-1\}.\]
For each $m\in M$ the connection $\nab$ restricts to give a holomorphic connection
$\nab_{m}$ on a trivial bundle over the space
\[\C_m=\C\setminus\{u_0,\cdots,u_{n-1}\}.\]
Dubrovin showed that these connections $\nab_m$ vary isomonodromically. We briefly explain this condition.

Choose a point $m\in M$ and a loop $\gamma_{m}$ in
$\C_{m}$ based at some point $z\in \C_m$. Let $H$ be the space of flat sections of $\nab_m$ near $z\in \C_m$. Monodromy around the loop $\gamma_m$ defines a linear transformation $\alpha_m\in \Aut(H)$. For points $m'\in M$ in a small neighbourhood of $m$ the connection $\nab$ allows us to identify $H$ with the space of flat sections of the connection $\nab_{m'}$ near $z\in \C_{m'}$. Moreover we can continuously deform the loop $\gamma_{m}$ to give a loop $\gamma_{m'}$ in $\C_{m'}$ based at $z$, and hence obtain a transformation $\alpha_{m'}\in \Aut(H)$. 
The isomonodromy condition is the statement that 
the transformations $\alpha_{m'}$ of $H$ obtained in this way are constant.

\subsection{Quantum cohomology of $\PP^2$}

Let us now consider the Frobenius manifold defined by the quantum cohomology of $\PP^2$. It is known that the subset of tame points of the resulting Frobenius manifold is non-empty, so we can apply Dubrovin's result Theorem \ref{dub} to obtain a tame Frobenius manifold structure on a dense open subset $M=\conft_3(\C)\setminus B$ where
$\conft_3(\C)$
is the universal cover. 
Let $\nab=\nab^{(\frac{1}{2})}$ be the second structure connection with parameter $s=1/2$.
The following result of Dubrovin's computes its monodromy.

\begin{thm}[Dubrovin]
\label{hann}
There is a point $m\in M$ with canonical co-ordinates $(u_0,u_1,u_2)$ the three roots of unity. Let $\gamma_0,\gamma_1,\gamma_2$ be the following loops $\gamma_i$ in $\C_m$ based at $0\in \C$
\[\includegraphics{paths.eps}\]
There is  a triple $(\phi_0,\phi_1,\phi_2)$ of flat sections of the connection $\nab_{m}$ in a neighbourhood of $0\in \C$,
such that the monodromy transformations $P_i$ corresponding to 
the loops $\gamma_i$ act on the triple $(\phi_0,\phi_1,\phi_2)$ by the matrices
\[\bigmat{1}{3}{-3}{0}{1}{0}{0}{0}{1},\quad \bigmat{1}{0}{0}{-3}{1}{3}{0}{0}{1},
\quad\bigmat{1}{0}{0}{0}{1}{0}{3}{-3}{1}.\]
Moreover this triple $(\phi_0,\phi_1,\phi_2)$ is unique up to multiplication by an overall scalar factor.
\end{thm}

\begin{pf}
The existence of a triple of flat sections with the above monodromy properties follows from general work of Dubrovin on monodromy of twisted period maps \cite[Lemma 4.10, 4.12]{Du5}, together with Dubrovin's computation of the Stokes matrix of the quantum cohomology of $\PP^2$
\cite[Example 4.4]{Du3}. Uniqueness is easily checked.
\end{pf}

The discriminant of the Frobenius manifold $M$ is the submanifold
\[\Delta=\{m\in M: u_i(m)=0\text{ for some }i\}.\]
Write $M^0=M\setminus\Delta$ for its complement and let $\Mt^0$ be its inverse image under the natural map
$\conft_3(\C^*)\to \conft_3(\C)$.
Let us take the point $m\in M^0$ of Theorem \ref{hann} as a base-point, and choose a small simply-connected neighbourhood $m\in U\subset M^0$. For each point $m\in U$ we have a well-defined choice of loops $\gamma_i$ in $\C_m$ based at $0$ obtained by deforming the loops $\gamma_i$ of Theorem \ref{hann}.

The group $G$ is a subgroup of the fundamental group $\pi_1(\conf_3(\C^*))$ and therefore acts by covering transformations on $\conft_3(\C^*)$.  Thus to each element $g\in G$ we can associate a corresponding open subset \[U_g=g(U)\cap \Mt^0\subset\Mt^0.\]
Using the connection $\nab$ we can continue the triple of sections of Theorem \ref{hann} to obtain a standard triple of flat sections of $\nab_{m'}$ in a neighbourhood of $0\in \C$ for all $m'\in \Mt^0$. This is well-defined despite the fact that $\Mt^0$ may not be simply-connected because of the isomonodromy property and the uniqueness statement in Theorem \ref{hann}. In particular, for each $g\in G$ we obtain a standard triple of flat sections of $\nab_{m'}$ near $0\in \C$ for all points $m'\in U_g$.  Taking monodromy of the connection $\nab_{m'}$ around the loops $\gamma_i$ with respect to this standard triple gives three matrices $P_0(g), P_1(g), P_2(g)$.

To calculate these matrices we use the isomonodromy property. For example, consider a path in $\Mt^0$ from $m$ to $\tau_1(m)$. If we move the loops $\gamma_i$ continuously with the $u_i$ then at the point $\tau_1(m)$ we will obtain the following basis $(\gamma'_0,\gamma'_1,\gamma'_2)$ of $\pi_1(\C_m,0)$.
\[\includegraphics{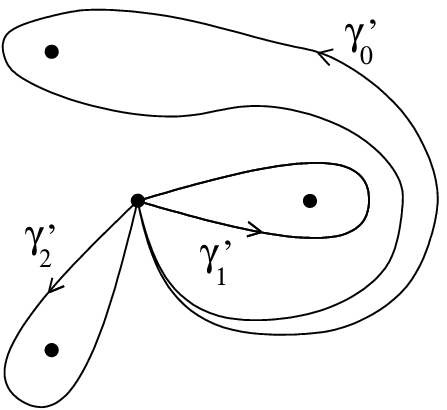}\]
Clearly
\[\gamma'_0=\gamma_0^{-1}\gamma_1 \gamma_0,\quad \gamma'_1=\gamma_0,\quad \gamma'_2=\gamma_2.\]
By the isomonodromy property, the monodromy of the standard triple of sections at $\tau_1(m)$ around the loops $\gamma'_i$ will be the same as the monodromy of the triple of sections at $m$ around the loops $\gamma_i$ and is therefore described by the matrices $P_i(e)$ of Theorem \ref{hann}. But the matrices $P_i(\tau_1)$ describe the monodromy of the same sections around the loops $\gamma_i$. Arguing in this way we see that the matrices $P_i(g)$ have the transformation properties
\begin{eqnarray*}
P_0(\tau_1 g)=P_1(g), \quad &P_1(\tau_1 g)=P_1(g) P_0(g) P_1(g)^{-1}, \quad &P_2(\tau_1 g)=P_2(g), \\
P_0(rg)=P_2(g), \quad &P_1(rg)=P_0(g), \quad &P_2(rg)=P_0(g).
\end{eqnarray*}
These are exactly the same transformation properties satisfied by the linear maps $\phi_{S_i(g)}$. Since the matrices $P_i(e)$ of Theorem \ref{hann} coincide with the matrices of $\phi_{S_i(e)}$ with respect to the basis $\{S_0(e),S_1(e),S_2(e)\}$ we obtain Theorem \ref{two}.

\subsection{Twisted period map}
There is an embedding
$M^0\hookrightarrow W$
obtained by sending a point $m$ to $(m,0)$. Pulling back the flat connection $\nab$ we obtain a flat connection on the tangent bundle $\T_{M^0}$ which we also denote $\nab$. We define flat co-ordinates $\WW_i$ whose gradients with respect to the flat metric on $M$ are the flat sections $(\phi_0,\phi_1,\phi_2)$ of Theorem \ref{hann}. Putting them together gives a holomorphic map
\[\WW\colon \Mt^0\lra \C^3\]
uniquely defined up to scalar multiples.
There is a free action of $\C$ on $\conft_3(\C^*)$ lifting the $\C^*$ action which simultaneously rescales the co-ordinates $(u_0,u_1,u_2)$ on $\conf_3(\C^*)$.
Let
\[\AA^2=\{(z_0,z_1,z_2)\in \C^3: z_0+z_1+z_2=i\}\]
be the affine space defined in the introduction. Then

\begin{prop}
There is a unique scalar multiple of the map $\WW$ which descends to give a local isomorphism
\[\WW\colon \Mt^0/\C\lra \AA^2.\]
\end{prop}

\begin{pf}
First we show that the only possible linear relation between the solutions $\phi_i$ of Theorem \ref{hann} is
\[\phi_0+\phi_1+\phi_2=0.\]
Indeed, any such relation must be monodromy invariant, and $(1,1,1)$ is the unique vector (up to multiples) preserved by the three given matrices.

Secondly we show that this relation does indeed hold. Otherwise $(\phi_0,\phi_1,\phi_2)$ define a basis of solutions and the map $\WW$ is a local isomorphism.
Let \[E=\big.\sum_i u_i \frac{\partial}{\partial u_i}\]
be the Euler vector field.
Dubrovin showed that all components $\WW_i$ of $\WW$ satisfy
\[\Lie_E(\WW_i)=\text{constant}.\] 
We cannot have $\Lie_E(\WW)=0$ since this would contradict the statement that $\WW$ is a local isomorphism. Thus there is a two-dimensional subspace of solutions satisfying $\Lie_E(\WW_i)=0$. But this subspace would have to be monodromy invariant, and there are no such two-dimensional subspaces. This gives a contradiction, so the relation holds, and rescaling we can assume that
\[\WW_0+\WW_1+\WW_2=i.\]
Now it follows that the only two-dimensional, monodromy invariant subspace of solutions is that generated by the $\phi_i$, so that $\Lie_E(\WW)=0$ and the result follows.
\end{pf}

\bigskip

\noindent School of Mathematics,
University of Sheffield,
Hicks Building, Hounsfield Road, Sheffield, S3 7RH, UK.

\smallskip

\noindent email: {\tt t.bridgeland@sheffield.ac.uk}


\begin{thebibliography}{101}

\bibitem{AGM} P. Aspinwall, B. Greene and D. Morrison, Measuring small distances in N=2 sigma models, Nuclear Phys. B 420 (1994), no. 1--2, 184--242.

\bibitem{Bo}
A. Bondal,
Representations of associative algebras and coherent
sheaves. (Russian) Izv. Akad. Nauk SSSR Ser. Mat.  53  (1989),
no. 1, 25--44;  translation in  Math. USSR-Izv.  34  (1990),
no. 1, 23--42

\bibitem{BP}
A. Bondal and A. Polishchuk,
Homological properties of associative algebras: the method of helices.
(Russian)  Izv. Ross. Akad. Nauk Ser. Mat.  57  (1993),  no. 2, 3--50;
translation in  Russian Acad. Sci. Izv. Math.  42  (1994), 
no. 2, 219--260.

\bibitem{BB} S. Brenner and M. Butler,
Generalizations of the Bernstein-Gelfand-Ponomarev reflection
functors. Representation theory, II (Proc. Second Internat. Conf.,
Carleton Univ., Ottawa, Ont., 1979),  pp. 103--169,
Lecture Notes in Math., 832, Springer, Berlin-New York, 1980.

\bibitem{B} {T. Bridgeland,} Stability conditions on triangulated categories, math.AG/0212237,
to appear in Ann. of Maths.

\bibitem{Br} T. Bridgeland, Stability conditions on K3 surfaces, math.AG/0307164.

\bibitem{B2} {T. Bridgeland,} T-structures on some local Calabi-Yau varieties,
J. of Algebra, 289 (2005) 453--483.

\bibitem{Br4} T. Bridgeland, Stability conditions and Kleinian singularities, math.AG/0508257.

\bibitem{BKR} T. Bridgeland, A. King and M. Reid, The McKay
correspondence as an equivalence of derived categories,
J. Amer. Math. Soc. {\bf 14} (2001), no. 3, 535-554.

\bibitem{CV} S. Cecotti and C. Vafa, On classification of $N=2$ supersymmetric
theories,  Comm. Math. Phys. 158 (1993), no. 3, 569--644.

\bibitem{DG} D.-E. Diaconescu and J. Gomis, Fractional branes and
boundary states in orbifold theories,  J. High Energy Phys. 2000, no. 10, Paper 1, 44 pp.

\bibitem{Do2} M. Douglas, B. Fiol and C. R{\"o}melsberger, The
spectrum of BPS branes on a noncompact Calabi-Yau, preprint
hep-th/0003263.

\bibitem{Do4} M. Douglas, Dirichlet branes,
homological mirror symmetry, and stability, Proceedings of the International Congress of Mathematicians, Vol. III (Beijing, 2002), 395--408, Higher Ed. Press, Beijing, 2002.


\bibitem{Du1} B. Dubrovin, Geometry of $2$D topological field theories. Integrable systems and quantum groups (Montecatini Terme, 1993), 120--348, Lecture Notes in Math., 1620, Springer, Berlin, 1996.

\bibitem{Du2} B. Dubrovin, Geometry and analytic theory of Frobenius manifolds. Proceedings of the International Congress of Mathematicians, Vol. II (Berlin, 1998)  Doc. Math. 1998. 

\bibitem{Du3} B. Dubrovin, Painlev{\'e} transcendents in two-dimensional topological field theory. The Painlev{\'e} property, 287--412, CRM Ser. Math. Phys., Springer, New York, 1999. 

\bibitem{Du4} B. Dubrovin and M. Mazzocco, Monodromy of certain Painlev{\'e}-VI transcendents and reflection groups. Invent. Math. 141 (2000), no. 1, 55--147.

\bibitem{Du5} B. Dubrovin, On almost duality for Frobenius manifolds. Geometry, topology, and mathematical physics, 75--132, Amer. Math. Soc. Transl. Ser. 2, 212, Amer. Math. Soc., Providence, RI, 2004.

\bibitem{FHHI} B. Feng, A. Hanany, Y. He and A. Iqbal, Quiver
theories, soliton spectra and Picard-Lefschetz transformations,  J. High Energy Phys. 2003, no. 2, 056, 33 pp.

\bibitem{GR} A. Gorodentsev and A. Rudakov,
Exceptional vector bundles on projective spaces,
Duke Math. J. 54 (1987), no. 1, 115--130.

\bibitem{Gu} D. Guzzetti,
Stokes matrices and monodromy of the quantum cohomology of projective spaces,
Comm. Math. Phys. 207 (1999), no. 2, 341--383.


\bibitem{HRS} D. Happel, I. Reiten and S. Smal{\o}, Tilting in abelian
categories and quasitilted algebras. Mem. Amer. Math. Soc. 120 (1996),
no. 575.


\bibitem{He}  C. Hertling, Frobenius manifolds and moduli spaces for singularities. Cambridge Tracts in Mathematics, 151. Cambridge University Press, Cambridge, 2002.

\bibitem{KP} R. Kent IV and D. Peifer, A geometric and algebraic description of annular braid groups, Internat. J. Algebra Comput. 12  (2002) 85--97.

\bibitem{Ko} M. Kontsevich, Homological algebra of mirror symmetry. Proceedings of the International Congress of Mathematicians, Vol. 1, 2 (Z{\"u}rich, 1994), 120--139, Birkhäuser, Basel, 1995.

\bibitem{Ma} Yu. Manin, Frobenius manifolds, quantum cohomology, and moduli spaces. American Mathematical Society Colloquium Publications, 47. American Mathematical Society, Providence, RI, 1999.

\bibitem{ST} P. Seidel and R. Thomas, Braid group actions on derived
categories of coherent sheaves, Duke Math. J. 108 (2001), no. 1,
37--108.

\bibitem{Ta} S. Tanab{\'{e}}, 
Invariant of the hypergeometric group associated to the quantum cohomology of the projective space,  Bull. Sci. Math. 128 (2004), no. 10, 811--827.


\bibitem{Za} E. Zaslow, Solitons and helices: the search for a math-physics bridge, 
Comm. Math. Phys. 175 (1996), no. 2, 337--375.


\end{thebibliography}
\end{document}